\documentclass{amsart}
\usepackage{amsthm,amssymb,amsmath}
\usepackage{epsfig}
\usepackage{color}

\newtheorem{prelem}{{\bf Proposition}}

\newtheorem{theorem}{Theorem}
\newtheorem{corollary}[theorem]{Corollary}

\newtheorem{proposition}[theorem]{Proposition}

\theoremstyle{definition}

\theoremstyle{remark}

\input{tcilatex}

\begin{document}
\title{roman domination and mycieleki's structure in graphs}
\author{ Adel P. Kazemi \vspace{4mm}\\
Department of Mathematics\\
University of Mohaghegh Ardabili\\
P. O. Box 5619911367, Ardabil, Iran\\
adelpkazemi@yahoo.com\vspace{3mm} \\
}
\date{}
\maketitle

\begin{abstract}
For a graph $G=(V,E)$, a function $f:V\rightarrow \{0,1,2\}$ is
called Roman dominating function (RDF) if for any vertex $v$ with
$f(v)=0$, there is at least one vertex $w$ in its neighborhood with
$f(w)=2$. The weight of an RDF $f$ of $G$ is the value
$f(V)=\sum_{v\in V}f(v)$. The minimum weight of an RDF of $G$ is its
Roman domination number and denoted by $\gamma_ R(G)$. In this
paper, we first show that $\gamma _{R}(G)+1\leq \gamma _{R}(\mu
(G))\leq \gamma _{R}(G)+2$, where $\mu (G)$ is the Mycielekian graph
of $G$, and then characterize the graphs achieving equality in these
bounds. Then for any positive integer $m$, we compute the Roman
domination number of the $m$-Mycieleskian $\mu _{m}(G)$ of a special
Roman graph $G$ in terms on $\gamma_R(G)$. Finally we present
several graphs to illustrate the discussed graphs.
\end{abstract}



\section{introduction and primary results }

The research of the domination in graphs has been an evergreen of
the graph theory. Its basic concept is the dominating set and the
domination number. The recent book $Fundamentals$ $of$ $Domination$
$in$ $Graphs$ \cite{HHS1} lists, in an appendix, many varieties of
dominating sets that have been studied. It appears that none of
those listed are the same as Roman dominating sets. Thus, Roman
domination appears to be a new variety of both historical and
mathematical interest.

A subset $S\subseteq V(G)$ is a $dominating$ $set$, briefly DS, in
$G$, if every vertex in $V(G)-S$ has a neighbor in $S$. The minimum
number of vertices of a DS in a graph $G$ is called the $domination$
$number$ of $G$ and denoted by $\gamma (G)$.

Let $f\colon V\rightarrow \{0,1,2\}$ be a function and let $%
(V_{0},V_{1},V_{2})$ be the ordered partition of $V$ induced by $f$, where $%
V_{i}=\{v\in V\mid f(v)=i\}$ and $\mid V_{i}\mid =n_{i}$, for
$i=0,1,2$. We notice that there is an obvious one-to-one
correspondence between $f$ and the
ordered partition $(V_{0},V_{1},V_{2})$ of $V$. Therefore, one can write $%
f=(V_{0},V_{1},V_{2})$. Function $f=(V_{0},V_{1},V_{2})$ is a
\emph{Roman
dominating function}, abbreviated RDF, if $V_{0}\subseteq N_{G}(V_{2})$. If $%
W_{2}\subseteq V_{2}$ and $W_{1}\subseteq V_{1}$, then we say
$W_{1}\cup W_{2}$ \emph{defends} $W_{1}\cup N_{G}[W_{2}]$. For
simplicity in notation,
instead of saying that $\{v\}$ defends $\{w\}$, we say $v$ \emph{defends} $w$%
. The \emph{weight} of $f$ is the value $f(V)=\sum_{v\in V}f(v)=2n_{2}+n_{1}$%
. The \emph{Roman domination number} $\gamma _{R}(G)$ is the minimum
weight of an RDF of $G$, and we say a function
$f=(V_{0},V_{1},V_{2})$ is a $\gamma _{R}$-\emph{function} if it is
an RDF and $f(V)=\gamma _{R}(G)$.

Stated in other words, a Roman dominating function is a coloring of
the vertices of a graph with the colors $\{0,1,2\}$ such that every
vertex colored $0$ is adjacent to at least one vertex colored $2$.
The definition of a Roman dominating function is given implicitly in
\cite{RR00} and \cite {St99}. The idea is that colors $1$ and $2$
represent either one or two Roman legions stationed at a given
location (vertex $v$). A nearby location (an adjacent vertex $u$) is
considered to be unsecured if no legions are stationed there (i.e.
$f(u)=0$). An unsecured location $(u)$ can be secured by sending a
legion to $u$ from an adjacent location $(v)$. But Emperor
Constantine the Great, in the fourth century A.D., decreed that a
legion cannot be sent from a location $v$ if doing so leaves that
location unsecured (i.e. if $f(v)=1$). Thus, two legions must be
stationed at a location $(f(v)=2)$ before one of the legions can be
sent to an adjacent location. More details about domination number
have given in many papers. For example reader can
see~\cite{CDHH04,FKKS09,He02,LKLP08,RR00, St99}.

As we will see, the generalized Mycieleskian graphs, which are also
called \emph{cones over graphs} \cite{Ta01}, are natural
generalization of Mycieleski graphs. If
$V(G)=V^{0}=\{v_{1}^{0},v_{2}^{0},\ldots ,v_{n}^{0}\}$ and
$E(G)=E_{0}$, then for any integer $m\geq 1$ the
$m$-\emph{Mycieleskian} $\mu _{m}(G)$ of $G$ is the graph with
vertex set $V^{0}\cup V^{1}\cup V^{2}\cup \cdots \cup V^{m}\cup
\{u\}$, where $V^{i}=\{v_{j}^{i}\mid
v_{j}^{0}\in V^{0}\}$ is the $i$-th distinct copy of $V^{0}$, for $%
i=1,2,\ldots ,m$, and edge set $E_{0}\cup \left( \cup
_{i=0}^{m-1}\{v_{j}^{i}v_{j^{\prime }}^{i+1}\mid
v_{j}^{0}v_{j^{\prime
}}^{0}\in E_{0}\}\right) \cup \{v_{j}^{m}u\mid v_{j}^{m}\in V^{m}\}.$ The $1$%
-\emph{Mycieleskian} $\mu _{1}(G)$ of $G$ is the same $Mycieleskian$
of $G$, and denoted simply by $\mu (G)$ or $M(G)$. Interested
readers may refer to \cite{BaRa08,LWLG,Ta01,We01} to know more about
the Mycieleskian graphs.

As stated in many references, for example in \cite{HHS1}, the \emph{%
Cartesian product} $G\,\Box \,H$ of two graphs $G$ and $H$\ is the
graph with vertex set $V(G)\times V(H)$ where two vertices $(u_{1},v_{1})$ and $%
(u_{2},v_{2})$ are adjacent if and only if either $u_{1}=u_{2}$ and $%
v_{1}v_{2}\in E(H)$ or $v_{1}=v_{2}$ and $u_{1}u_{2}\in E(G)$.

The notation we use is as follows. Let $G$ be a simple graph with $vertex$ $%
set$ $V=V(G)$ and $edge$ $set$ $E=E(G)$. The $order$ $\mid V\mid $
and $size$ $\mid E\mid $\ of $G$ are respectively denoted by
$n=n(G)$ and $m=m(G)$. For every vertex $v\in V$, the $open$ $neighborhood$ $N_{G}(v)$ is the set $%
\{u\in V\mid uv\in E\}$ and the $closed$ $neighborhood$ of $v$ is the set $%
N_{G}[v]=N_{G}(v)\cup \{v\}$. Also for a subset $X\subseteq V(G)$, the \emph{open neighborhood} of $X$ is $%
N_{G}(X)=\cup _{v\in X}N_{G}(v)$ and its \emph{closed neighborhood} is $%
N_{G}[X]=N_{G}(X)\cup X$. The $degree$ of a vertex $v\in V$ is $%
deg(v)=\mid N(v)\mid $. The $minimum$ and $maximum$ $degree$ of a
graph $G$ are denoted by $\delta =\delta (G)$ and $\Delta =\Delta
(G)$, respectively. If every vertex of $G$ has degree $k$, then $G$
is said to be $k$-$regular$. We write $K_{n}$, $C_{n}$ and $P_{n}$,
respectively, for the $complete$ $graph$, $cycle$ and $path$ of
order $n$ and $K_{n_{1},...,n_{p}}$ for the $complete$ $p$-$partite$
$graph$.

Let $v\in S\subseteq V$. A vertex $u$ is called a \emph{private
neighbor} of $v$ with respect to $S$, or simply an $S$-$pn$ of $v$,
if $u\in N[v]-N[S-\{v\}]$. The set $pn(v;S)=N[v]-N[S-\{v\}]$\ of all
$S$-pn's of $v$ is called the \emph{private neighborhood set} of $v$
with respect to $S$.
Also an $S$-pn of $v$ is an \emph{external private neighbor} or \emph{%
external} (denoted by $S$-$epn$ of $v$) if it is a vertex of $V-S$.
We also call the set $epn(v;S)=N(v)-N[S-\{v\}]$\ of all $S$-epn's of $v$, the \emph{%
external private neighborhood set} of $v$ with respect to $S$. To
see this definitions refer to \cite{BaRa08,HHS1}. Obviously if
$f=(V_{0},V_{1},V_{2})$ is a $\gamma _{R}$-function, then for each
$v\in V_{2}$, $epn(v;V_{2})\neq \emptyset$ (we notice that for each
vertex $v\in V_2$, $epn(v;V_{2})\subseteq V_0$ and so
$epn(v;V_{2})\neq \emptyset$ if and only if $epn(v;V_{2})\cap V_0
\neq \emptyset$).

Cockayne et al. in~\cite{CDHH04} have shown that for any graph $G$ of order $%
n$ and maximum degree $\Delta $, $2n/(\Delta +1)\leq \gamma
_{R}(G)$, and for the classes of paths $P_{n}$ and cycles $C_{n}$,
$\gamma _{R}(P_{n})=\gamma _{R}(C_{n})=\lceil 2n/3\rceil$.
Furthermore, they have shown that for any graph $G$, $\gamma (G)\leq
\gamma _{R}(G)\leq 2\gamma (G)$, where the lower bound is achieved
only by $G=\overline{K_{n}}$, the empty graph with $n$ vertices. A
graph $G$
is called a \emph{Roman graph} if $\gamma _{R}(G)=2\gamma (G)$~\cite{CDHH04}%
. For example, the complete multipartite graph $K_{m_{1},\ldots
,m_{n}}$ is Roman if and only if $2\notin \{m_{1},\ldots ,m_{n}\}$. As shown in \cite{CDHH04}, an equivalent condition for $G$ to be a Roman graph is that $%
G $ has a $\gamma _{R}$-function $f=(V_{0},V_{1},V_{2})$ with $%
V_{1}=\emptyset$.

We now introduce two new concepts. A Roman graph $G$ with $\gamma
_{R}$-function $f=(V_{0},\emptyset,V_{2})$ we call a \emph{special
Roman graph} if the induced subgraph $G[V_{2}]$ has no isolated
vertex, and its $\gamma
_{R}$-function $f=(V_{0},\emptyset,V_{2})$ we call a \emph{special} $\gamma _{R}$-%
\emph{function}.

In this paper, we first show that $\gamma _{R}(G)+1\leq \gamma
_{R}(\mu (G))\leq \gamma _{R}(G)+2$ and characterize the graphs
achieving equality in these bounds. Then for any positive integer
$m$, we compute the $\gamma_R(\mu _m(G))$ of a special Roman graph
$G$ in terms on $\gamma_R(G)$. Finally we present several graphs to
illustrate the discussed graphs.

In this entire of paper we assume that the induced subgraph by
$V_{1}$ is an empty subgraph if $V_{1}\neq \emptyset$.

We first present the Roman domination number of some known graphs.

\begin{prelem}
\label{Com.Multi.Rom} \emph{(\textbf{Cockayne et al. \cite{CDHH04}
2004})} If
$m_{1}\leq m_{2}\leq \cdots \leq m_{n}$ are positive integers and $%
G $ is the complete $n$-partite graph $K_{m_{1},\ldots ,m_{n}}$,
then\textrm{
\[
\gamma _{R}(G)=\left\{
\begin{array}{ll}
m_{1}+1 & \mbox{if }\;1\leq m_{1}\leq 2, \\
4 & \mbox{otherwise.}
\end{array}
\right.
\]
}
\end{prelem}

\begin{prelem}
\label{PtxKn.Rom} Let $t\geq 1$ and $n\geq 3$ be integers. If $G$ is
the cartesian product graph $P_t\square K_n$, then
\textrm{
\[
\gamma _{R}(G)=\left\{
\begin{array}{ll}
6\lfloor t/4\rfloor +2r & \mbox{if }\;n=3, \mbox{ }t \equiv r \mbox{ } (\mbox{mod }4) \mbox{ and }0\leq r\leq 2, \\
6\lfloor t/4\rfloor +2r-1 & \mbox{if }\;n=3, \mbox{ }t \equiv 3 \mbox{ } (\mbox{mod }4), \\
2t & \mbox{otherwise.}
\end{array}
\right.
\]
}
\end{prelem}

\begin{proof}
Let $V(P_t\square K_n)=\{1,2,...,t\}\times \{1,2,...,n\}$. We first
suppose that $n=3$. Let $A=\{(4\ell +1,1),(4\ell +3,2)\mid 0\leq
\ell \leq \lfloor t/4\rfloor-1 \}$ and $B=\{(4\ell +2,3),(4\ell
+4,3)\mid 0\leq \ell \leq \lfloor t/4\rfloor-1 \}$. Easily it can be
seen that the given Roman dominating functions have minimum weight.

\textbf{Case i}. $t \equiv 0 \mbox{ (mod }4)$.\\
Let $f_0=(W_0,W_1,W_2)$, where $W_2=A$, $W_1=B$ and $W_0=V-(W_1\cup
W_2)$.

\textbf{Case ii}. $t \equiv 1 \mbox{ (mod }4)$.\\
Let $f_1=(W_0^{'},W_1^{'},W_2^{'})$, where $W_2^{'}=A\cup
\{(t,1)\}$, $W_1^{'}=B$ and $W_0^{'}=V-(W_1^{'}\cup W_2^{'})$.

\textbf{Case iii}. $t \equiv 2 \mbox{ (mod }4)$.\\
Let $f_2=(W_0^{''},W_1^{''},W_2^{''})$, where $W_2^{''}=A\cup
\{(t-1,1),(t,1)\}$, $W_1^{''}=B$ and $W_0^{''}=V-(W_1^{''}\cup
W_2^{''})$.

\textbf{Case iv}. $t \equiv 3 \mbox{ (mod }4)$.\\
Let $f_3=(W_0^{'''},W_1^{'''},W_2^{'''})$, where $W_2^{'''}=A\cup
\{(t-2,1),(t,2)\}$, $W_1^{''}=B\cup \{(t-1,3)\}$ and
$W_0^{'''}=V-(W_1^{'''}\cup W_2^{'''})$.\\
Now let $n\geq 4$. Easily it can be seen that the wight of every RDF
for $P_t\square K_n$ on the every copy of $K_n$ is at least 2. Thus
$\gamma _{R}(P_t\square K_n)\geq 2t$. Now since
$f=(W_0,\emptyset,W_2)$ is an RDF with weight $2t$, when
$W_2=\{(\ell ,1)\mid 1\leq \ell \leq t \}$ and $W_0=V-W_2$, then
$\gamma _{R}(P_t\square K_n)= 2t$.
\end{proof}

\begin{prelem}
\label{CtxKn.Rom} Let $t\geq 1$ and $n\geq 3$ be integers. If $G$ is
the cartesian product graph $C_t\square K_n$, then\textrm{
\[
\gamma _{R}(G)=\left\{
\begin{array}{ll}
6\lfloor t/4\rfloor +2r & \mbox{if }\;n=3, \mbox{ }t \equiv r \mbox{ } (\mbox{mod }4) \mbox{ and }0\leq r\leq 1, \\
6\lfloor t/4\rfloor +2r-1 & \mbox{if }\;n=3, \mbox{ }t \equiv r \mbox{ } (\mbox{mod }4) \mbox{ and }2\leq r\leq 3, \\
2t & \mbox{otherwise.}
\end{array}
\right.
\]
}
\end{prelem}

\begin{proof}
Let $V(C_t\square K_n)=\{1,2,...,t\}\times \{1,2,...,n\}$. We first
suppose that $n=3$. Let $A=\{(4\ell +1,1),(4\ell +3,2)\mid 0\leq
\ell \leq \lfloor t/4\rfloor-1 \}$ and $B=\{(4\ell +2,3),(4\ell
+4,3)\mid 0\leq \ell \leq \lfloor t/4\rfloor-1 \}$. Easily it can be
seen that the given Roman dominating functions have minimum weight.

\textbf{Case i}. $t \equiv 0 \mbox{ (mod }4)$.\\
Let $f_0=(W_0,W_1,W_2)$, where $W_2=A$, $W_1=B$ and $W_0=V-(W_1\cup
W_2)$.

\textbf{Case ii}. $t \equiv 1 \mbox{ (mod }4)$.\\
Let $f_1=(W_0^{'},W_1^{'},W_2^{'})$, where $W_2^{'}=A\cup
\{(t,1)\}$, $W_1^{'}=B$ and $W_0^{'}=V-(W_1^{'}\cup W_2^{'})$.

\textbf{Case iii}. $t \equiv 2 \mbox{ (mod }4)$.\\
Let $f_2=(W_0^{''},W_1^{''},W_2^{''})$, where $W_2^{''}=A\cup
\{(t-1,3)\}$, $W_1^{''}=(B-\{(t-2,3)\})\cup \{(t-2,1),(t,2)\}$ and
$W_0^{''}=V-(W_1^{''}\cup W_2^{''})$.

\textbf{Case iv}. $t \equiv 3 \mbox{ (mod }4)$.\\
Let $f_3=(W_0^{'''},W_1^{'''},W_2^{'''})$, where $W_2^{'''}=A\cup
\{(t-2,1),(t,2)\}$, $W_1^{''}=B\cup \{(t-1,3)\}$ and
$W_0^{'''}=V-(W_1^{'''}\cup W_2^{'''})$.\\
Similar to the proof of Proposition \ref{PtxKn.Rom}, we can proof
that $\gamma _{R}(C_t\square K_n)= 2t$, when $n\geq 4$.
\end{proof}

Similarly, the following two propositions can be proved and easily
can be verified that the given graphs in them are special Roman
graphs.

\begin{prelem}
\label{PtxKn,..,n.Rom} 
 If $t\geq 2$ and $4\leq n_1\leq n_2\leq ... \leq n_p$ are integers, then $\gamma _{R}(P_t\square
 K_{n_1,...,n_p})=4t$.
\end{prelem}

\begin{prelem}
\label{PtxK1,n.Rom} 
 If $t\geq 1$ and $n\geq 2$ are integers, then $\gamma _{R}(P_t\square
 K_{1,n})=2t$.
\end{prelem}

\section{Main Results}

First we state our main theorem.

\begin{theorem}
\label{ULbounds} For each graph $G$, $\gamma_{R}(G)+1\leq
\gamma_{R}(\mu (G))\leq \gamma_{R}(G)+2$.
\end{theorem}

\noindent

\begin{proof}
Let $V(G)=V^{0}$, and $V(\mu (G))=V^{0}\cup V^{1}\cup \{u\}$. Let $%
f=(W_{0},W_{1},W_{2})$\ be a $\gamma _{R}$-function of $G$. Since $%
g=(W_{0}\cup V^{1},W_{1},W_{2}\cup \{u\})$ is an RDF of $\mu (G)$, $%
\gamma _{R}(\mu (G))\leq \gamma _{R}(G)+2$. We now show that $\gamma
_{R}(\mu
(G))\geq \gamma _{R}(G)+1$. Let $g=(W_{0},W_{1},W_{2})$\ be an $\gamma _{R}$%
-function of $\mu (G)$. We continue our discussion in the following
two cases.

\textbf{Case 1.} $u\in W_{1}\cup W_{2}$. Let

\[
\begin{array}{lcl}
W_{2}^{\prime } & = & (W_{2}-(\{u\}\cup (W_{2}\cap V^{1})))\cup
\{v_{j}^{0}\mid v_{j}^{1}\in W_{2}\}, \\
W_{1}^{\prime } & = & W_{1}-\{v_{j}^{0}\mid v_{j}^{1}\in W_{2}\},
\\
W_{0}^{\prime } & = & V(\mu (G))-(W_{1}^{\prime }\cup W_{2}^{\prime }).
\end{array}
\]

\noindent Then, the function $g^{\prime }=(W_{0}^{\prime },W_{1}^{\prime
},W_{2}^{\prime })$ is an RDF for $G$, and hence
\[
\gamma _{R}(G)\leq g^{\prime }(V(G))=2|W_{2}^{\prime }|+|W_{1}^{\prime
}|\leq 2|W_{2}|+|W_{1}|-1\leq \gamma _{R}(\mu (G))-1,
\]
and so $\gamma _{R}(\mu (G))\geq \gamma _{R}(G)+1$, as desired.

\textbf{Case 2.} $u\in W_{0}$.

\noindent Then $W_{2}\cap V^{1}\neq \emptyset$. If $W_{1}\cap
V^{1}\neq \emptyset$, then the inequality is easily seen to be true. So let $%
W_{1}\subseteq V^{0}$. Consider $A=\{v_{j}^{0}|v_{j}^{1}\in W_{2}\}$. If $%
(W_{2}\cup W_{1})\cap A\neq \emptyset$, then similar to Case 1, we
can find an RDF $g^{\prime }$ for $G$ such that $\gamma _{R}(G)\leq
g^{\prime }(V(G))\leq \gamma _{R}(\mu (G))-1$, and hence $\gamma
_{R}(\mu (G))\geq \gamma _{R}(G)+1$. Therefore,
we assume that $(W_{2}\cup W_{1})\cap A=\emptyset$ and $%
W_{2}=\{v_{s_{i}}^{0}|1\leq i\leq t\}\cup \{v_{j_{\ell}}^{1}|1\leq
\ell \leq m\}$\ such that $(W_{1}\cup \{v_{s_{i}}^{0}|1\leq i\leq
t\})\cap \{v_{j_{\ell}}^{0}|1\leq \ell \leq m\}=\emptyset$.

Assume that $epn(v_{j_{k}}^{1};W_{2})\cap V^{0}=\emptyset$, for some $%
1\leq k\leq m$. Then, $k$ is unique and
$epn(v_{j_{k}}^{1};W_{2})=\{u\}$. Let
\[
\begin{array}{lcl}
W_{2}^{\prime } & = & \{v_{s_{i}}^{0}|1\leq i\leq t\}\cup
\{v_{j_{\ell}}^{0}|1\leq \ell \leq m,\mbox{ and }\ell \neq k\}, \\
W_{1}^{\prime } & = & W_{1}, \\
W_{0}^{\prime } & = & V(G)-(W_{1}^{\prime }\cup W_{2}^{\prime }).
\end{array}
\]
Then the function $g^{\prime }=(W_{0}^{\prime },W_{1}^{\prime
},W_{2}^{\prime })$ is an RDF for $G$ such that $\gamma _{R}(G)\leq
\gamma _{R}(\mu (G))-2$, which implies that $\gamma _{R}(\mu
(G))\geq \gamma _{R}(G)+1$, as desired. Hence we may assume that for
each $1\leq \ell\leq m,\ $ $epn(v_{j_{\ell}}^{1};W_{2})\cap
V^{0}\neq \emptyset$.

Let $\alpha _{j_{\ell}}^{0}\in epn(v_{j_{\ell}}^{1};W_{2})\cap V^{0},$ for every $\ell$, $%
1\leq \ell \leq m$. Clearly $\{\alpha _{j_{l}}^{\ell}| 1\leq \ell
\leq m\}\cap (W_{0}\cup W_{1})=\emptyset$. Therefore $\{\alpha
_{j_{\ell}}^{1}| 1\leq \ell \leq m\}\subseteq W_{2}$. Further $m\geq
2$, and $\{\alpha _{j_{\ell}}^{1}| 1\leq \ell \leq
m\}=\{v_{j_{\ell}}^{1}| 1\leq \ell \leq m\}$. Also for any $\ell$,
$1\leq \ell \leq m$, $\mid epn(v_{j_{\ell}}^{1};W_{2})\cap V^{0}\mid
=1$. Let $\alpha
_{j_{1}}^{1}=v_{j_{2}}^{1}$. We now add $v_{j_{2}}^{0}$\ and $v_{j_{2}}^{1}$%
\ to $W_{2}$\ and $W_{1}$, respectively, and delete $v_{j_{1}}^{1}$\ and $%
v_{j_{2}}^{1}$\ of $W_{2}$. If necessary, we also add $u$ to $W_{1}$. Then
we obtain $g^{\prime }=(W_{0}^{\prime },W_{1}^{\prime },W_{2}^{\prime }) $
as a new RDF for $\mu (G)$. If $m\geq 3$, then $g^{\prime }(u)=0$, and hence
$g^{\prime }(V(\mu (G))\leq g(V(\mu (G)))-1=\gamma_{R}(\mu (G))-1, $\ a
contradiction. Finally let $m=2$ and choose $W_{2}^{\prime \prime
}=W_{2}^{\prime }$, $W_{1}^{\prime \prime }=W_{1}^{\prime
}-\{u,v_{j_{2}}^{1}\}$, and $W_{0}^{\prime \prime }=V(G)-(W_{1}^{\prime
\prime }\cup W_{2}^{\prime \prime })$. The function $g^{\prime \prime
}=(W_{0}^{\prime \prime },W_{1}^{\prime \prime },W_{2}^{\prime \prime })$\
is an RDF for $G$ with weight $\gamma_{R}(\mu (G))-2$. Therefore, $%
\gamma_{R}(G)\leq g^{\prime \prime }(V(G))\leq \gamma _{R}(\mu (G))-2$, as
desired.
\end{proof}

\medskip Our next aim is to characterize for which graphs $G$ the Roman
domination number of $\mu (G)$\ is $\gamma _{R}(G)+1$ or $\gamma _{R}(G)+2$.

\begin{theorem}
\label{+1} For any special Roman graph $G$, $\gamma_{R}(\mu
(G))=\gamma _{R}(G)+1.$
\end{theorem}

\noindent

\begin{proof}
By Theorem \ref{ULbounds}, $\gamma _{R}(\mu (G))\geq \gamma _{R}(G)+1$. Let $%
f=(V_{0},\emptyset,V_{2})$\ be a special $\gamma _{R}$-function for
$G$. By choosing $W_{2}=V_{2}$, $W_{1}=\{u\}$, and $W_{0}=V_{0}\cup
V^{1}$, we see that the function $g=(W_{0},W_{1},W_{2})$ is an RDF
for $\mu (G)$\ with weight $\gamma _{R}(G)+1$, which implies that
$\gamma _{R}(\mu (G))=\gamma _{R}(G)+1$.
\end{proof}

\medskip In the next theorem we show that the converse of Theorem \ref{+1} is also
true.

\begin{theorem}
\label{+2} If $G$\ is not a Roman graph, then $\gamma_{R}(\mu
(G))=\gamma _{R}(G)+2.$
\end{theorem}

\noindent

\begin{proof}
In the contrary, let $g=(W_{0},W_{1},W_{2})$\ be a $\gamma
_{R}$-function for $\mu (G)$ with weight $\gamma _{R}(G)+1$, by
Theorem 1. We assume that if $\mid W_{1}\mid \geq 1$, then the
induced subgraph $\mu (G)[W_{1}]$\ is isomorphic to the empty graph
$\overline{K_{b}}$, where $b=\mid W_{1}\mid$. In the next three
cases, we show that $u\notin W_{0}\cup W_{1}\cup W_{2}$ which
completes our proof.

\textbf{Case 1.} $u\in W_{2}$.

\noindent Then $W_{1}\subseteq V^{0}$. Let
\[
\begin{array}{lcl}
V_{2} & = & (W_{2}-\{u\})\cup \{v_{j}^{0}|v_{j}^{1}\in W_{2}\}, \\
V_{1} & = & W_{1}-\{v_{j}^{0}|v_{j}^{1}\in W_{2}\}, \\
V_{0} & = & V^{0}-(V_{1}\cup V_{2}).
\end{array}
\]

\noindent Then the function $f=(V_{0},V_{1},V_{2})$\ is an RDF for $G$ with
at most weight $\gamma _{R}(G)-1$, a contradiction.

\textbf{Case 2. }$u\in W_{1}$.

\noindent Then $W_{2}\subseteq V^{0}$. Since the induced subgraph $\mu
(G)[W_{1}]$\ is an empty graph, we have $W_{1}-\{u\}\subseteq V^{0}.$ Since $%
v_{j}^{0}\in W_{1}$ implies $v_{j}^{1}\in W_{1}$, thus $W_{1}=\{u\}$. Let $%
V_{2}=W_{2}$, $V_{1}=\emptyset$, and $V_{0}=V^{0}-V_{2}$. Then the
function $f=(V_{0},\emptyset,V_{2})$\ is a $\gamma _{R}$-function
for $G$ and hence $G$ is a Roman graph that is a contradiction.

\textbf{Case 3.} $u\in W_{0}$.

\noindent Then $| W_{2}\cap V^{1}| \geq 1$. Let $v_{1}^{1}\in
W_{2}\cap V^{1} $. Then we may assume that $v_{1}^{0}\in W_{0}\cup
W_{1}$. Because if $v_{1}^{0}\in W_{2}$, then with considering
\[
\begin{array}{lcl}
V_{2} & = & (W_{2}\cap V^{0})\cup \{v_{j}^{0}|v_{j}^{1}\in W_{2}\}, \\
V_{1} & = & (W_{1}\cap V^{0})-\{v_{j}^{0}|v_{j}^{1}\in W_{2}\},
\\
V_{0} & = & V^{0}-(V_{1}\cup V_{2}),
\end{array}
\]
the function $f=(V_{0},V_{1},V_{2})$\  is an RDF for $G$ with at
most weight $\gamma_{R}(\mu (G))-2=\gamma_{R}(G)-1$ that is a
contradiction. We now continue our discussion on the following two
subcases.

\textbf{Subcase 3.i. }$v_{1}^{0}\in W_{1}$.

\noindent Let $A=\{v_{j}^{1}|v_{j}^{0}\in W_{1}\}$. Since
$\{i|v_{i}^{0}\in W_{1}$ and $v_{i}^{1}\in W_{0}\}=\emptyset$, and
$g$ is a $\gamma _{R}$ -function for $\mu (G)$, thus $|A\cap
W_{2}|\leq 1$. Since $v_{1}^{1}\in W_{2}\cap A$,
then $|A\cap W_{2}|=1$. Easily we see that if $v_{i}^{1}\in W_{1}$, then $%
v_{i}^{0}\in W_{1}\cup W_{2}$. Let $t$ be the number of $i$ s which
$v_{i}^{0},v_{i}^{1}\in W_{1}$ and let $\ell$ be the number of $i$ s
which $v_{i}^{0}\in W_{2}$ and $v_{i}^{1}\in W_{1}$. If $t+ \ell
\geq1$, then we can get an RDF with weight $\gamma _R(G)-1$ for $G$.
Therefore, let $\ell=t=0$. Thus $W_{1}=\{v_{1}^{0}\}$ and
$epn(v_{1}^{1};W_{2})=\{u\}$. Then with considering
\[
\begin{array}{lcl}
V_{2} & = & (W_{2}\cap V^{0})\cup \{v_{j}^{0}|v_{j}^{1}\in
W_{2}\}-\{v_{1}^{0}\}, \\
V_{1} & = & W_{1}, \\
V_{0} & = & V^{0}-(V_{1}\cup V_{2})
\end{array}
\]
the function $f=(V_{0},V_{1},V_{2})$\ is an RDF for $G$ such that
$f(V(G))\leq \gamma_{R}(\mu (G))-2=\gamma_{R}(G)-1$, that is a
contradiction.

\textbf{Subcase 3.ii. }$v_{1}^{0}\in W_{0}$.

\noindent We recall $v_{1}^{1}\in W_{2}$ and $u\in W_{0}$. By
Subcase 3.i and above discussion we may assume that if $v_{i}^{1}\in
W_{2}$, then $v_{i}^{0}\in W_{0}$. We also know that
if $v_{i}^{0}\in W_{2}$, then $v_{i}^{1}\notin W_{2}$ . The assumption $%
v_{1}^{0}\in W_{0}$\ concludes that $v_{1}^{0}$\ is defended by a vertex $%
\alpha $ of $W_{2}$. Suppose $\alpha \in V^{0}$ and let $\alpha
=v_{2}^{0}$. If $epn(v_{1}^{1};W_{2})\cap V^{0}=\emptyset$, then,
with deleting at least $v_{1}^{1}$\ from $W_{2}$, we can define a
function $f=(V_{0},V_{1},V_{2})$\
with weight at most $\gamma _{R}(\mu (G))-2=\gamma _{R}(G)-1$\ such that $%
V_{1}\cup V_{2}$\ defends all vertices of $G$. Thus let $%
epn(v_{1}^{1};W_{2})\cap V^{0}\neq \emptyset$, and for some $t\geq
3$ let
$epn(v_{1}^{1};W_{2})\cap V^{0}=\{v_{i}^{0}|3\leq i\leq t\}$. Then $%
\{v_{i}^{1}|3\leq i\leq t\}\subseteq W_{1}\cup W_{2}$. If $g(\cup
_{i=3}^{t}v_{i}^{1})\geq 2$, then improving $g$ makes a function
$g^{\prime } $ for $\mu (G)$\ with weight at most $\gamma _{R}(\mu
(G))-1$. Then $epn(v_{1}^{1};W_{2})\cap V^{0}=\{v_{3}^{0}\}$. In
this case, we may find an RDF $g^{\prime }$ for $\mu (G)$\ such that either $%
g^{\prime }(V(\mu (G)))\leq \gamma _{R}(\mu (G))-1$, which is a
contradiction, or $g^{\prime }(V(\mu (G)))=\gamma _{R}(\mu (G))$ and $%
g^{\prime }(u)=1$, that is impossible by Case 2. Finally we assume that $%
W_{2}\cap V^{0}$\ does not defend $v_{1}^{0}$ and let $\alpha
=v_{2}^{1}$. Then $epn(v_{2}^{1};W_{2})\cap V^{0}\neq \emptyset$.
Also we have $epn(v_{1}^{1};W_{2})\cap V^{0}\neq \emptyset$. If
$epn(v_{2}^{1};W_{2})\cap V^{0}=\{v_1^0\}$, then with choosing
\[
\begin{array}{lcl}
V_{2}^{'} = (W_{2}\cap V^{0})\cup \{v_{j}^{0}|v_{j}^{1}\in
W_{2}\}-\{v_{1}^{1}\}, \\
V_{1}^{'} = (W_{1}\cap V^{0}) - \{v_{j}^{0}\in W_1 \mid
v_{j}^{1}\in W_{2}\}, \\
V_{0}^{'} = V(G)-(V_{1}^{'}\cup V_{2}^{'}),
\end{array}
\]
the function $g^{'}=(V_{0}^{'},V_{1}^{'},V_{2}^{'})$\ is an RDF for
$G$ such that $g^{'}(V(G))\leq g(V(G))-2 = \gamma_{R}(G)-1$, that is
a contradiction. Thus let $\mid epn(v_{2}^{1};W_{2})\cap V^{0}\mid
\geq 2$. Then $v_j^{1}\in W_1 \cup W_2$, when $v_j^{0}\in
epn(v_2^{1};W_2)\cap V^{0}$. Then with choosing
\[
\begin{array}{lcl}
W_{2}^{'} =  (W_{2}-\{v_1^{1},v_2^{1}\})\cup
\{v_1^{0},v_2^{0}\}, \\
W_{1}^{'} =  W_{1} - \{v_{j}^{1} \mid
v_{j}^{0}\in epn(v_2^{1};W_2) \}, \\
W_{0}^{'} =  V(\mu (G))-(W_{1}^{'}\cup W_{2}^{'}),
\end{array}
\]
the function $g^{'}=(W_{0}^{'},W_{1}^{'},W_{2}^{'})$ is an RDF for
$\mu (G)$ such that
\[
g^{'}(V(\mu (G)))=2\mid W_2^{'} \mid + \mid W_1^{'} \mid \leq 2\mid
W_2 \mid + \mid W_1 \mid -1= g(V(\mu (G)))-1 = \gamma_{R}(\mu
(G))-1,
\]
that is a contradiction.
\end{proof}

The proof of the next theorem is similar to Theorem \ref{+2} and
hence is removed for brevity.

\begin{theorem}
\label{mu+2} For any $\gamma _{R}$-function $f=(V_{0},V_{1},V_{2})$
of a graph $G,$ if the induced subgraph $G[V_{2}]$ has an isolated
vertex, then $\gamma _{R}(\mu (G))=\gamma _{R}(G)+2$.
\end{theorem}

The next two theorems are immediate results of Theorems \ref{+1} and
\ref{mu+2}.

\begin{theorem}
\label{mu} Let $G$ be any graph. Then $\gamma _{R}(\mu (G))=\gamma
_{R}(G)+1$ if and only if $G$ is a special Roman graph.
\end{theorem}

\begin{theorem}
\label{mu2} Let $G$ be any graph. Then $\gamma _{R}(\mu (G))=\gamma
_{R}(G)+2$ if and only if $G$ is not a special Roman graph.
\end{theorem}

For any integer $m\geq 2$, we now compute the Roman domination number of the
$m$-Mycieleskian $\mu_{m}(G)$ of a special Roman graph $G$.

\begin{theorem}
\label{mu.m} Let $m$ be a positive integer. If $G$ is a special
Roman graph, then
\[
\gamma_{R}(\mu_{m}(G))=\left\{
\begin{array}{lcl}
2\lceil m/4\rceil\gamma_{R}(G)+2 & \mbox{if } m\equiv 0 \mbox{ } (\mbox{mod }4), \\
(2\lceil m/4\rceil -1)\gamma_{R}(G)+1 & \mbox{if } m\equiv 1
 \mbox{ } (\mbox{mod }4), \\
2\lceil m/4\rceil\gamma_{R}(G) & \mbox{if } m\equiv 2  \mbox{ } (\mbox{mod }4), \\
2\lceil m/4\rceil\gamma_{R}(G) +1 & \mbox{if } m\equiv 3  \mbox{ }
(\mbox{mod }4).
\end{array}
\right.
\]
\end{theorem}

\noindent

\begin{proof}
Let $f=(V_{0},\emptyset,V_{2})$ be a special $\gamma
_{R}$-function\ for $G$. Then $\gamma _{R}(G)\geq 4$, and for each $v\in V_{2}$, $%
\mid epn(v;V_{2})\mid \geq 2$, and hence $\mid V(G)\mid \geq 3\mid
V_2 \mid = 3\gamma _{R}(G)/2$. We will make an RDF
$g=(W_{0},W_{1},W_{2})$ for $\mu _{m}(G)$ with minimum weight.

Let $g=(W_{0},W_{1},W_{2})$ be an RDF for $\mu _{m}(G)$ with minimum
weight. First we show that the number of vertices of $W_{1}\cup
W_{2}$ that dominate all vertices of four consecutive rows among the
rows
 $V^{0}$, $V^{1}$, $\dots$ , $V^{m}$ is $\gamma _{R}(G)$.
 Let $2\leq i\leq m-2$, and let $S$ be a subset of $W_{1}\cup W_{2}$
 that dominates all vertices of  $V^{i-1}$. Then $S\subseteq V^{i-2}\cup V^{i}$.
 Without loss of generality, let $S\subseteq V^{i}$. Since the projection of $S$ on $%
V^{0}$ dominates all vertices of $V^{0}$ and $G$ is Roman, thus
$\mid S\mid \geq \gamma _{R}(G)/2$. We notice that $V_2^i$ with
cardinal $\gamma _{R}(G)/2$ dominates all vertices of $V^{i-1}$,
where $V_2^i=\{v^i\mid v \in V_2\}$. Thus $S$ can be $V_2^{i}$. But
$S$ dominates no vertices of $V^i-S$. Let $t=\mid V^{i}-S\mid $.
Hence $t \geq 2\mid S \mid =\gamma_ R(G)$. Because $\mid V(G)\mid
=\mid V^i\mid \geq (3/2)\gamma_R(G)$ and $\mid S \mid
=\gamma_R(G)/2$. Now for dominating $V^{i}-S$, we choose a subset
$S^{\prime }$ of $V^{i+1}$ such that it dominates $V^i$ (and hence
$V^i-S$) and $V^{i+2}$. Similar to above discussions, we have $\mid
S^{\prime } \mid \geq \gamma_R(G)/2$, and we may assume that $\mid
S^{\prime } \mid = \gamma_R(G)/2$. Thus $S \cup S^{\prime }$ is a
dominating set of the vertices of $V^{i-1}\cup V^{i}\cup V^{i+1}\cup
V^{i+2}$. Now we discuss on $r$, where $m\equiv r \mbox{ (mod }4)$.

Let $\lfloor m/4\rfloor =k$. We recall that
$f=(V_{0},\emptyset,V_{2})$\ is a special $\gamma _{R}$-function\
for $G$ and $V_{2}^{i}$ denotes the set $\{v^{i}\in V^{i}\mid v\in
V_{2}\}$.

For $m\equiv 0 \mbox{ ( mod }4)$, let

\[
\begin{array}{llll}
W_{2} & = & \cup _{t=0}^{k-1}(V_2^{4t+1}\cup V_2^{4t+2})\cup \{u\},
&
\mbox {where }k\geq 1, \\
W_{1} & = & \emptyset, &  \\
W_{0} & = & V(\mu _{m}(G))-(W_{1}\cup W_{2}). &
\end{array}
\]

For $m\equiv 1 \mbox{ ( mod }4)$, let

\[
\begin{array}{llll}
W_{2} & = & \cup _{t=0}^{k-1}(V_2^{4t+1}\cup V_2^{4t+2})\cup V_{2}^{m}, & %
\mbox {where }k\geq 0, \\
W_{1} & = & \{u\}, &  \\
W_{0} & = & V(\mu _{m}(G))-(W_{1}\cup W_{2}). &
\end{array}
\]

For $m\equiv 2 \mbox{ ( mod }4)$, let

\[
\begin{array}{llll}
W_{2} & = & \cup _{t=0}^{k}(V_2^{4t+1}\cup V_2^{4t+2}), \mbox {  where }k\geq 0, \\
W_{1} & = & \emptyset , &  \\
W_{0} & = & V(\mu _{m}(G))-(W_{1}\cup W_{2}). &
\end{array}
\]

For $m\equiv 3 \mbox{ ( mod }4)$, let

\[
\begin{array}{llll}
W_{2} & = & \cup _{t=0}^{k}(V_2^{4t+1}\cup V_2^{4t+2}), \mbox {  where }k\geq 0, \\
W_{1} & = & \{u\} , &  \\
W_{0} & = & V(\mu _{m}(G))-(W_{1}\cup W_{2}).  &
\end{array}
\]
Then the function $g=(W_{0},W_{1},W_{2})$ is an RDF with minimum
weight for $\mu _{m}(G)$\ such that
\[
g((\mu_{m}(G)))=\left\{
\begin{array}{lcl}
2\lceil m/4\rceil\gamma_{R}(G)+2 & \mbox{if } m\equiv 0 \mbox{ ( mod }4), \\
(2\lceil m/4\rceil -1)\gamma_{R}(G)+1 & \mbox{if } m\equiv 1
\mbox{ ( mod }4), \\
2\lceil m/4\rceil\gamma_{R}(G) & \mbox{if } m\equiv 2 \mbox{ ( mod }4), \\
2\lceil m/4\rceil\gamma_{R}(G) +1 & \mbox{if } m\equiv 3 \mbox{ (
mod }4),
\end{array}
\right.
\]
as desired.
\end{proof}

One can verify that the given graphs in the next three propositions
are special Roman graphs, and by Theorem \ref{mu.m} and Propositions
\ref{Com.Multi.Rom}, \ref{PtxKn.Rom}, \ref{CtxKn.Rom},
\ref{PtxKn,..,n.Rom} and \ref{PtxK1,n.Rom}, they are proved. In the
next propositions, it is assumed that $m\geq 2$.

\begin{proposition}
\label{mu.m.multi} Let $3\leq m_{1}\leq m_{2}\leq \cdots \leq m_{n}$
be integers. If $G=K_{m_{1},\ldots ,m_{n}}$, then
\[
\gamma _{R}(\mu _{m}(G))=\left\{
\begin{array}{lcl}
2m+2 & \mbox{if } m\equiv 0 \mbox{ }( \mbox{mod }4), \\
2m+3 & \mbox{if } m\equiv 2 \mbox{ }( \mbox{mod }4),\\
2m+4 & \mbox{otherwise.}
\end{array}
\right.
\]
\end{proposition}

\begin{proposition}
\label{mu.m.cart.1} If $G\in \{P_{t}\square K_{n}|t\geq 2,\mbox{ }
n\geq 4\}\cup \{C_{t}\square K_{n}|t\geq 3,\mbox{ } n\geq 4\}\cup
\{P_{t}\square K_{1,n}|t\geq 2,\mbox{ } n\geq 2\}$, then
\[
\gamma _{R}(\mu _{m}(G))=\left\{
\begin{array}{lcl}
mt+2 & \mbox{if } m\equiv 0 \mbox{ }( \mbox{mod }4), \\
(m+2)t & \mbox{if } m\equiv 2 \mbox{ }( \mbox{mod }4),\\
(m+1)t+1 & \mbox{otherwise.}
\end{array}
\right.
\]
\end{proposition}

\begin{proposition}
\label{mu.m.cart.2} Let $t\geq 2$ and $4\leq n_1\leq ... \leq n_p$
be integers. If $G=P_{t}\square K_{n_1,...,n_p}$, then
\[
\gamma _{R}(\mu _{m}(G))=\left\{
\begin{array}{lcl}
2mt+2 & \mbox{if } m\equiv 0 \mbox{ }( \mbox{mod }4), \\
2(m+2)t & \mbox{if } m\equiv 2 \mbox{ }( \mbox{mod }4),\\
2(m+1)t+1 & \mbox{otherwise.}
\end{array}
\right.
\]
\end{proposition}

There are some graphs that are not special Roman graph. The complete graphs $%
K_{n}$, paths $P_{n}$, stars $K_{1,n}$, all for $n\geq 1$, cycles $C_{n}$
for $n\geq 3$, and complete multipartite graphs $K_{2,m_{2},\ldots ,m_{n}}$
for $2\leq m_{2}\leq \cdots \leq m_{n}$ are not special Roman graphs and
their Roman domination number are respectively $2$, $\lceil 2n/3\rceil $, $2$%
, $\lceil 2n/3\rceil $, and\ $3$. The next proposition gives another non
special Roman graph.

\begin{proposition}
\label{petersen} The Petersen graph $G(5)$ is not a special Roman
graph, and $\gamma _{R}(G(5))=6$.
\end{proposition}

\noindent

\begin{proof}
Let $V(G(5))=\{i|1\leq i\leq 10\}$ and
\[
E(G(5))=\{(i,i+1)|1\leq i\leq 9\}\cup
\{(6,10),(1,5),(1,9),(2,7),(3,10),(4,8)\}.
\]
Since $G(5)$\ is $3$-regular and $G(5)$\ has 10 vertices, then $\gamma
_{R}(G(5))\geq 6$. Let $V_{2}=\{1,8,10\}$, $V_{1}=\emptyset$, and $%
V_{0}=V-V_{2}$. Since the function $f=(V_{0},\emptyset,V_{2})$\ is
an RDF with weight $6$, then $\gamma _{R}(G(5))=6$.

Finally we prove $G(5)$\ is not a special Roman graph. By the given RDF in
the previous paragraph, $G(5)$ is a Roman graph. Now let $%
f=(V_{0},\emptyset,V_{2})$\ be an arbitrary $\gamma _{R}$-function\ for $%
G(5)$. We know $G(5)$\ is a non-complete $3$-partite graph with three parts $%
X=\{1,4,7,10\}$, $Y=\{3,6,8\}$, and $Z=\{2,5,9\}$. Assume that $a$
and $b$ are two adjacent vertices of $V_{2}$. Since $a$ and $b$
defend together six vertices, and there is no other vertex $c$ that
dominates all of the four remained vertices, then $f(V)\geq 7$, a
contradiction. Hence $G(5)$ is not\ a special Roman graph.
\end{proof}

\begin{corollary}
\label{examp+2} Let G be $K_{n}$, $P_{n}$, $K_{1,n}$ for $n\geq 1$, $C_{n}$ for $n\geq 3$,\ $%
K_{2,m_{2},\ldots,m_{n}}$ for $2\leq m_{2} \leq m_{3} \leq \cdots
\leq m_{n}$ or $G(5)$. Then $\gamma_{R}(\mu (G))=\gamma_{R}(G)+2.$
\end{corollary}


\section {acknowledgements}

We thanks the referee for his/her helpful suggestions.


\begin{thebibliography}{99}
\bibitem{BaRa08}  R. Balakrishnan and S. Francis Raj, Connectivity of the
Mycielskian of a graph, \textit{Discrete Math}. \textbf{308} (2008),
2607-2610.

\bibitem{CDHH04}  E. J. Cockayne, P. A. Dreyer Jr., S. M. Hedetniemi, S. T.
Hedetniemi, Roman domination in graphs, \textit{Discrete Math.}
\textbf{278} (2004) 11-22.

\bibitem{FKKS09}  O. Favaron, H. Karami, R. Khoeilar, and S. M.
Sheikholeslami, On the Roman domination number of a graph,
\textit{Discrete Math.} \textbf{309} (2009), 3447-3451.

\bibitem{HHS1}  T.W. Haynes, S.T. Hedetniemi, and P.J. Slater. Fundamentals
of Domination in Graphs. Marcel Dekker, New York, 1998.

\bibitem{He02}  M. A. Henning, A Characterization of Roman trees,
\textit{Discussiones Mathematicae Graph Theory} \textbf{22} (2002),
325-334.

\bibitem{LKLP08}  M. Liedloff, T. Kloks, J. Liu, and S. L. Peng, Efficient
algorithms for Roman domination on some classes of graphs,
\textit{Discrete Appl. Math.} \textbf{156} (2008), 3400-3415.

\bibitem{LWLG}  W. Lin, J. Wu, P. C. B. Lam, and G. Gu, Several parameters
of generalized Mycielskians, \textit{Discrete Appl. Math.}
\textbf{154} (2006) 1173-1182.

\bibitem{RR00}  C. S. ReVelle, K.E. Rosing, Defendens imperium romanum: a
classical problem in military strategy, \textit{Amer. Math. Monthly}
\textbf{107} (7) (2000) 585-594.

\bibitem{St99}  I. Stewart, Defend the Roman Empire!, \textit{Sci. Amer.} \textbf{281} (6)
(1999) 136-139.

\bibitem{Ta01}  C. Tardif, Fractional chromatic numbers of cones over
graphs, \textit{J. Graph Theory} \textbf{38} (2001) 87-94.

\bibitem{We01}  D. B. West, Introduction to graph theory, (2nd edition),
\textit{Prentice Hall} USA (2001).
\end{thebibliography}
\end{document}